\documentclass{svjour3A}                     
\smartqed  
\usepackage{graphicx}
\usepackage{amsmath}
\usepackage{amsfonts,amssymb}

\newcommand{\zr}{\ltimes}
\newcommand{\Real}{\mathbb{R}}
\newcommand{\Co}{\mathbb{C}}

\newcommand{\g}{\mathfrak{g}}
\newcommand{\h}{\mathfrak{h}}

\newcommand{\so}{\mathfrak{so}}

\def\sl{\mathfrak{sl}}
\newcommand{\gl}{\mathfrak{gl}}

\def\sp{\mathfrak{sp}}

\newcommand{\R}{\mathcal{R}}

\newcommand{\id}{\mathop\text{\rm id}\nolimits}

\newcommand{\pr}{\mathop\text{\rm pr}\nolimits}

\def\g{\mathfrak{g}}
\def\h{\mathfrak{h}}
\def\so{\mathfrak{so}}

\def\sl{\mathfrak{sl}}
\def\sp{\mathfrak{sp}}
\def\gl{\mathfrak{gl}}

\def\pr{\mathop\text{\rm pr}\nolimits}

\def\id{\mathop\text{\rm id}\nolimits}

\def\GL{\text{\rm GL}}
\def\SL{\text{\rm SL}}
\def\SO{\text{\rm SO}}

\def\U{\text{\rm U}}
\def\Sp{\text{\rm Sp}}
\def\Spin{\text{\rm Spin}}

\def\SU{\text{\rm SU}}

\def\zr{\ltimes}

\def\Real{\mathbb{R}}
\def\Co{\mathbb{C}}
\def\g{\mathfrak{g}}
\def\h{\mathfrak{h}}

\def\so{\mathfrak{so}}

\def\sl{\mathfrak{sl}}
\def\gl{\mathfrak{gl}}

\def\sp{\mathfrak{sp}}

\def\R{\mathcal{R}}
\def\W{\mathcal{W}}

\def\H{\mathbb {H}}
\def\id{\text{\rm id}}

\def\pr{\text{\rm pr}}

\def\Mat{\mathop\text{\rm Mat}\nolimits}

\begin{document}

\title{Holonomy groups of pseudo-quaternionic-K\"ahlerian manifolds  of non-zero scalar curvature
} 

\titlerunning{ }        

\author{Natalia I. Bezvitnaya}

\authorrunning{ } 

\institute{N.I. Bezvitnaya \at
             Department of Mathematics and Statistics, Faculty of
Science, Masaryk University in Brno, Kotl\'a\v rsk\'a~2, 611~37
Brno, Czech Republic. \\
              \email{bezvitnaya@math.muni.cz}           
}

\date{ }

\maketitle

\centerline{\it Dedicated to Dmitri Vladimirovich Alekseevsky at
the occasion of his 70th birthday}

\vskip0.3cm

\begin{abstract}
 The holonomy group $G$ of a pseudo-quaternionic-K\"ahlerian manifold
of signature $(4r,4s)$ with non-zero scalar curvature is contained
in $\Sp(1)\cdot\Sp(r,s)$ and it contains $\Sp(1)$.  It is proved
that either $G$ is irreducible, or $s=r$ and $G$ preserves an
 isotropic subspace of dimension $4r$, in the last case, there are only two possibilities
 for the connected component of the identity of such $G$. This gives the classification
 of possible connected holonomy
 groups of pseudo-quaternionic-K\"ahlerian manifolds
of non-zero scalar curvature.

\end{abstract}

\section{Introduction and the results}

Quaternionic-K\"ahlerian geometry is of increasing interest both
in mathematics and mathematical physics, see e.g.
\cite{AC05,March,Salamon,Swann}. A pseudo-quaternionic-K\"ahlerian
manifold  is a pseudo-Riemannian manifold $(M,g)$ of signature
$(4r,4s)$, $r+s>1$ together with a parallel quaternionic structure
$Q\subset\so(TM)$, i.e. three-dimensional linear Lie algebra $Q$
with a basis $I_{1},I_{2},I_{3}$ which satisfies the relations
$I_{1}^{2}=I_{2}^{2}=I_3^2=-\id$, $I_{3}=I_{1}I_{2}=-I_{2}I_{1}$.
The holonomy group $G$ of such manifold is contained in
$\Sp(1)\cdot\Sp(r,s)=\Sp(1)\times\Sp(r,s)/\mathbb{Z}_2$.
Conversely, any  pseudo-Riemannian manifold with such holonomy
group is pseudo-quaternionic-K\"ahlerian. Denote by $G^0$ the
restricted holonomy group of $(M,g)$, i.e. the connected component
of identity of $G$.

The classification of connected holonomy groups of Riemannian
manifolds is well known and it has a lot of applications both in
geometry and physics, see e.g. \cite{Ber,Besse,Bryant,Jo}. The
corresponding problem for pseudo-Riemannian manifolds of arbitrary
signature is solved only in some partial cases, see the recent
reviews \cite{ESI,IRMA}. The difficulty appears if the holonomy
group preserves a degenerate subspace of the tangent space. Here
we show that the holonomy group $G\subset\Sp(1)\cdot\Sp(r,s)$ of a
pseudo-quaternionic-K\"ahlerian manifold with non-zero scalar
curvature is irreducible if $s\neq r$. If $s=r$, then $G$ may
preserve a degenerate subspace of the tangent space, in this case
there are only two possibilities for~$G$.

In \cite{AC05} it is proved  that the curvature tensor $R$ of a
pseudo-quaternionic-K\"ahlerian manifold of signature $(4r,4s)$
can be written as
\begin{equation}\label{RR0W}R(X,Y)=\nu R_0+{\cal W},\end{equation}
where $\nu=\frac{{\rm scal}}{4m(m+2)}$ ($m=4r+4s$) is the reduced
scalar curvature,
\begin{equation}\label{R0}R_0(X,Y)=\frac{1}{2}\sum_{\alpha=1}^3g(X,I_\alpha
Y)I_\alpha+\frac{1}{4}\left(X\wedge Y+\sum_{\alpha=1}^3I_\alpha
X\wedge I_\alpha Y\right),\end{equation} $X,Y\in TM$, is the
curvature tensor of the quaternionic projective space
$\mathbb{PH}^{r,s}$, and ${\cal W}$ is an algebraic curvature
tensor with zero Ricci tensor. It is proved that  ${\rm scal}\neq
0$ if and only if the holonomy  group $G$ contains $\Sp(1)$.
 In \cite{AC05} it is proved also that any
pseudo-quaternionic-K\"ahlerian manifold with non-zero scalar
curvature is locally indecomposable, i.e. it is not locally a
product of two pseudo-Riemannian manifolds of positive dimension,
or, equivalently, $G^0$ does not preserve any proper {\it
non-degenerate} vector subspace of the tangent space.

The tangent space to the manifold $(M,g)$ can be identified with
the pseudo-Euclidean space $\Real^{4r,4s}$ endowed with the
pseudo-Euclidean metric $\eta$ and an $\eta$-orthogonal
quaternionic structure $I_1,I_2,I_3$, or with the
pseudo-quaternionic-Hermitian space $\H^{r,s}$ endowed with a
pseudo-quaternionic-Hermitian metric $g$. Let $s=r$ and
$W\subset\H^{r,r}$ be an isotropic subspace of quaternionic
dimension $r$. Let $p_1,...,p_r,q_1,...,q_r$ be a basis of
$\H^{r,r}$ such that $p_1,...,p_r$ is a basis of $W$ and the only
non-zero values of the pseudo-quaternionic-Hermitian form $g$ on
$\H^{r,r}$ are $g(p_i,q_i)=g(q_i,p_i)=1$. The maximal subalgebra
$\sp(r,r)_W\subset\sp(r,r)$ preserving $W$ can be identified with
the following matrix Lie algebra:
$$\sp(r,r)_W=\left. \left\{\left(\begin{array}{cc}C&B\\0&-\bar
C^t\end{array}\right)\right|C,B\in\Mat(r,\H),\,\bar{B}^t=-B\right\},$$
where $\Mat(r,\H)$ denotes the space of $r\times r$ quaternionic
matrices. Note the following. Let $\mathbb{H}^{m}$ be an
m-dimensional quaternionic vector space and $e_{1},...,e_{m}$ a
basis of $\mathbb{H}^{m}$. We identify an element
$X\in\mathbb{H}^{m}$ with the column $(X_{t})$ of the left
coordinates of $X$ with respect to this basis,
$X=\sum_{t=1}^{m}X_{t}e_{t}$. Let
$f:\mathbb{H}^{m}\to\mathbb{H}^{m}$ be an $\mathbb{H}$-linear map.
Define the matrix $\Mat_{f}$ of $f$ by the relation
$fe_{l}=\sum_{t=1}^{m}(\Mat_{f})_{tl}e_{t}$. Now if
$X\in\mathbb{H}^{m}$, then $fX=(X^{t}\Mat_{f}^{t})^{t}$ and
because of the non-commutativity of the quaternionic numbers this
is not the same as $\Mat_{f}X$. Conversely, any $m\times m $
quaternionic matrix defines an $\mathbb{H}$-linear map
$f:\mathbb{H}^{m}\to\mathbb{H}^{m}$. Denote by
$\Sp(r,r)_W\subset\Sp(r,r)$ the maximal connected Lie subgroup
preserving $W$, i.e. the connected Lie subgroup corresponding to
the subalgebra $\sp(r,r)_W\subset\sp(r,r)$. Define the Lie
subalgebra
$$\h_0=\sp(1)\oplus\left.
\left\{\left(\begin{array}{cc}C&0\\0&-\bar
C^t\end{array}\right)\right|C\in\Mat(r,\H)\right\}\subset\sp(1)\oplus\sp(r,r)_W$$
and denote by $H_0$ the corresponding connected Lie subgroup of
$\Sp(1)\cdot\Sp(r,r)_W$.

We prove the following two statements.

\begin{theorem}\label{Main} Let $(M,g)$ be a pseudo-quaternionic-K\"ahlerian manifold
of non-zero scalar curvature and of signature $(4r,4s)$. If its
restricted holonomy group $G^0$ is not irreducible, then  $s=r$,
$G^0$ preserves an isotropic quaternionic subspace
$W\subset\H^{r,r}$ of quaternionic dimension $r$ and either
$G^0=H_0$, or $G^0=\Sp(1)\cdot\Sp(r,r)_W$.
\end{theorem}

\begin{proposition}\label{prop1} Any pseudo-quaternionic-K\"ahlerian manifold with
the restricted holonomy group  $H_0$ is locally symmetric, i.e.
its curvature tensor is parallel.
\end{proposition}

Note that if the manifold $(M,g)$ is not locally symmetric and its
restricted holonomy group $G^0$ is irreducible, then
$G^0=G=\Sp(1)\cdot\Sp(r,s)$, e.g. \cite{Bryant}.

Simply connected symmetric pseudo-quaternionic-K\"ahlerian
manifolds are classified in \cite{AC05}. Each such space $(M,g)$
may be represented as $M=F/K$, where $F$ is the connected group
generated by transvections and $K$ is the stabilizer of a fixed
point $o\in M$. The holonomy group of $(M,g)$ coincides with the
isotropy representation of $K$. These spaces are exhausted by
$$
\begin{array}{lllll}
\frac{\SU(p + 2, q)}{S(\U(2) \times \U(p, q))},& \frac{\SL(r +
1,\H)}{ S(\GL(1,\H) \times \GL(r,\H))},& \frac{\SO_0(p + 4,
q)}{\SO(4) \times \SO_0(p, q)},& \frac{ \SO^*(2l + 4)}{ \SO^*(4)
\times \SO^*(2l) },& \frac{\Sp(p + 1, q) }{\Sp(1) \times \Sp(p, q)
},\\  \frac{ E_{6(-78)}}{ \SU(2)\SU(6) },& \frac{ E_{6(2)}}{
\SU(2)\SU(6) },& \frac{ E_{6(2)}}{ \SU(2)\SU(2, 4) },& \frac{
E_{6(-14)}}{ \SU(2)\SU(2, 4) },& \frac{ E_{6(6)}}{ \Sp(1)\SL(3,\H)
},\\ \frac{ E_{6(-26)}}{ \Sp(1)\SL(3,\H)},& \frac{ E_{7(-133)}}{
\SU(2)\Spin(12) },& \frac{ E_{7(-5)}}{\SU(2)\Spin(12)},& \frac{
E_{7(-5)}}{ \SU(2)\Spin^0(4, 8)},& \frac{ E_{7(7)}}{
\SU(2)\SO^*(12) },\\ \frac{ E_{7(-25)}}{\SU(2)\SO^*(12)},&
\frac{E_{8(-248)}}{\SU(2)E_{7(133)}},& \frac{
E_{8(-24)}}{\SU(2)E_{7(133)}},&
\frac{E_{8(-24)}}{\SU(2)E_{7(-5)}},&
\frac{E_{8(8)}}{\SU(2)E_{7(-5)}},\\
\frac{F_{4(-52)}}{\Sp(1)\Sp(3)},& \frac{F_{4(4)}}{\Sp(1)\Sp(3)},&
\frac{F_{4(4)}}{\Sp(1)\Sp(1, 2)},& \frac{F_{4(-20)}}{\Sp(2)\Sp(1,
2) },& \frac{G_{2(-14)}}{\SO(4)},\quad \frac{
G_{2(2)}}{\SO(4)}.\end{array}$$ Note that the subgroup
$H_0\subset\Sp(1)\cdot\Sp(r,r)_W$ is the holonomy group of the
symmetric space $\frac{\SL(r + 1,\H)}{ S(\GL(1,\H) \times
\GL(r,\H))}$. This shows that {\it the holonomy groups of other
symmetric spaces are irreducible}.

We get the following corollaries.

\begin{corollary} Let $(M,g)$ be a  pseudo-quaternionic-K\"ahlerian manifold
of non-zero scalar curvature and of signature $(4r,4s)$. Then
either its restricted holonomy group coincides with
$\Sp(1)\cdot\Sp(r,s)$ or with $\Sp(1)\cdot\Sp(r,r)_W$, or $(M,g)$
is locally  symmetric.
\end{corollary}

\begin{corollary} Let $(M,g)$ be a  complete pseudo-quaternionic-K\"ahlerian manifold
of non-zero scalar curvature and of signature $(4r,4s)$. Then
either its restricted holonomy group coincides with
$\Sp(1)\cdot\Sp(r,s)$ or with $\Sp(1)\cdot\Sp(r,r)_W$, or $(M,g)$
is the factor space of a symmetric space obtained in \cite{AC05}
by a freely acting discrete group $\Gamma$.
\end{corollary}

\begin{corollary} Let $(M,g)$ be a simply connected complete pseudo-quaternionic-K\"ahlerian manifold
of non-zero scalar curvature and of signature $(4r,4s)$. If the
holonomy group $G$ of $(M,g)$ is irreducible, then either
$G=\Sp(1)\cdot\Sp(r,s)$, or $(M,g)$ is a symmetric space obtained
in \cite{AC05} different from $\frac{\SL(r + 1,\H)}{ S(\GL(1,\H)
\times \GL(r,\H))}$. If $G$ is not irreducible, then $s=r$, and
either $G=\Sp(1)\cdot\Sp(r,r)_W$, or $(M,g)$ is isometric to the
symmetric space $\frac{\SL(r + 1,\H)}{ S(\GL(1,\H) \times
\GL(r,\H))}$. \end{corollary}

 Theorem \ref{Main} gives
only the list of possible connected holonomy groups. To complete
the classification of all connected holonomy groups, one must show
that $\Sp(1)\cdot\Sp(r,r)_W$ may appear as the holonomy group of a
pseudo-quaternionic-K\"ahlerian manifold.

\section{Proof of Theorem \ref{Main}}

Let $(M,g)$ be a pseudo-quaternionic-K\"ahlerian manifold of
non-zero scalar curvature and of signature $(4r,4s)$. Let $m=r+s$.
Obviously, it is enough to prove the theorem in terms of the
holonomy algebra $\g\subset\sp(1)\oplus\sp(r,s)$ of $(M,g)$. In
\cite{AC05} it is proved that $\g$ contains $\sp(1)$ and it does
not preserve any proper non-degenerate subspace of
$\Real^{4r,4s}$.

For any subalgebra $\g\subset\sp(1)\oplus\sp(r,s)$ denote by
$\R(\g)$ the space of algebraic curvature tensors of type $\g$,
i.e. the space of linear maps from $\wedge^2\Real^{4r,4s}$ to $\g$
satisfying the first Bianchi identity
\begin{equation}\label{Bianchi1} R(X,Y)Z+R(Y,Z)X+R(Z,X)Y=0\end{equation} for all
$X,Y,Z\in\Real^{4r,4s}$. It is well-known that any $R\in\R(\g)$
satisfies
\begin{equation}\label{SYMR} \eta(R(X,Y)Z,U)=\eta(R(Z,U)X,Y)\end{equation} for all
$X,Y,Z,U\in\Real^{4r,4s}$.
 For
example, if $(M,g)$ is a pseudo-quaternionic-K\"ahlerian manifold,
$x\in M$, and $\g$ is the holonomy algebra of $(M,g)$ at the point
$x$, then identifying $T_xM$ with $\Real^{4r,4s}$, we get that the
value $R_x$ of the curvature tensor $R$ of $(M,g)$ at the point
$x$ belongs to $\R(\g)$. From the Ambrose-Singer Theorem
\cite{Besse} it follows that if $\g\subset\sp(1)\oplus\sp(r,s)$ is
the holonomy algebra of a pseudo-quaternionic-K\"ahlerian
manifold, then it is {a Berger algebra}, i.e. $\g$ is spanned by
the images of the algebraic curvature tensors $R\in\R(\g)$. We
will prove the theorem assuming that
$\g\subset\sp(1)\oplus\sp(r,s)$ is a Berger algebra.

Suppose that $\g$ preserves a degenerate subspace
$W\subset\Real^{4r,4s}$. Then $\g$ preserves the isotropic
subspace $W\cap W^\bot$, i.e. we may assume that $W$ is isotropic.
Since $\g$ contains $\sp(1)$, $W$ is a quaternionic subspace of
$\H^{r,s}$. Let $\dim_\H W=t$. Then $W\subset W^\bot$ and $\dim_\H
W^\bot=m-t>t$. Let $E$ be a quaternionic subspace of $W^\bot$
complementary to $W$, then the restriction of $g$ to $E$ is
non-degenerate, let $(r_0,s_0)$ be its signature. We have
$W\subset E^\bot$. Let $W_1\subset E^\bot$ be any isotropic
subspace complementary to $W$. Clearly, $\dim_\H W_1=\dim_\H W$.
Let $p_1,...,p_t,e_1,...,e_{r_0+s_0},q_1,...,q_t$ be a basis of
$\H^{r,s}$ such that $p_1,...,p_t\in W$, $e_1,...,e_{r_0+s_0}\in
E$, $q_1,...,q_t\in W_1$, and the only non-zero values of the
pseudo-quaternionic-Hermitian form $g$ are
$g(p_i,q_i)=g(q_i,p_i)=1$ (if $1\leq i \leq t$), $g(e_i,e_i)=-1$
(if $1\leq i \leq r_0$), $g(e_i,e_i)=1$ (if $r_0+1\leq i \leq
r_0+s_0$). Then the Lie algebra $\sp(r,s)$ can be identified with
the following matrix Lie algebra:
$$\sp(r,s)=\left.
\left\{\left(\begin{array}{ccc}C&-(E_{r_0,s_0}\bar{X})^t&B\\Y&A&X\\D&-(E_{r_0,s_0}\bar{Y})^t&-\bar{C}^t\end{array}\right)
\right|\begin{array}{cc}C\in\Mat(r,\H),& B,D\in S(r,\H),\\
A\in\sp(r_0,s_0),&X,Y\in\Mat(r_0+s_0,t,\H)\end{array} \right\},$$
where $\Mat(r,\H)$ denotes the space of $r\times r$ quaternionic
matrices, $$S(r,\H)=\{B\in\Mat(r,\H)|\,\ \bar{B}^t=-B\},$$
$\Mat(r_0+s_0,t,\H)$ denotes the space of $r_0+s_0\times t$
quaternionic matrices,
$E_{r_0,s_0}=\left(\begin{array}{cc}-E_{r_0}&0\\0&E_{s_0}\end{array}\right)$,
and $E_l$ denotes the identity $l\times l$ matrix.   For the
maximal subalgebra $\sp(r,s)_W\subset\sp(r,s)$ preserving $W$ we
get
$$\sp(r,s)_W=\left.
\left\{\left(\begin{array}{ccc}C&-(E_{r_0,s_0}\bar{X})^t&B\\0&A&X\\0&0&-\bar{C}^t\end{array}\right)
\right|\begin{array}{cc}C\in\Mat(r,\H),& B\in
S(r,\H),\\A\in\sp(r_0,s_0),& X\in\Mat(r_0+s_0,t,\H)\end{array}
\right\}.$$

We claim that if $r_0+s_0\neq 0$, then
$\R(\sp(1)\oplus\sp(r,s)_W)=\R(\sp(r,s)_W)$. Let
$R\in\R(\sp(1)\oplus\sp(r,s)_W)$. From \eqref{RR0W}
 it follows that $R=\nu R_0+\W,$ where $\nu\in\Real$, $R_0$ is given by
\eqref{R0} with $X,Y\in\Real^{4r,4s}$, and $\W\in\R(\sp(r,s))$.
Let $p\in W$, $X\in E$, and $Y,Z\in\Real^{4r,4s}$. Using
\eqref{SYMR}, we get
\begin{equation}\label{sym}\eta(R(p,X)Y,Z)=\eta(R(Y,Z)p,X)=0,\end{equation} since $R(Y,Z)p\in W$ and
$X\in E\subset W^\bot$. Consequently, $R(p,X)=0$. Let $X,Y\in E$,
then using the Bianchi identity, we get
$$R(X,Y)p=-R(Y,p)X-R(p,X)Y=0.$$ This shows that $\nu\,
\pr_W\circ R_0(X,Y)|_{W}=-\pr_W\circ\W(X,Y)|_{W}$. On the other
hand, $\pr_W\circ\W(X,Y)|_{W}\in\gl(W)=\gl(r,\H)$, whilst
$\pr_W\circ
R_0(X,Y)|_{W}=\frac{1}{2}\sum_{\alpha=1}^3\eta(X,I_\alpha
Y)I_\alpha|_{W}$. Consequently, $\nu \,\pr_W\circ
R_0(X,Y)|_{W}=-\pr_W\circ\W(X,Y)|_{W}=0$. This means that
$\nu\sum_{\alpha=1}^3\eta(X,I_\alpha Y)I_\alpha=0.$ Taking
$X=e_1$, $Y=I_1 e_1$, we get $\nu=0$. Thus,
$R=\W\in\R(\sp(r,s)_W)$. This shows that if $r_0+s_0\neq 0$, then
any Berger subalgebra of $\sp(1)\oplus\sp(r,s)_W$ is contained in
$\sp(r,s)_W$.

Let $r_0=s_0=0$. Then $s=t=r$. From \eqref{RR0W}
 it follows that \begin{equation}\label{Rsp1+}
\R(\sp(1)\oplus\sp(r,r))=\Real
R_0\oplus\R(\sp(r,r)),\end{equation} where $R_0$ is given by
\eqref{R0} with $X,Y\in\Real^{4r,4r}$. Consider the subalgebra
\begin{equation}\h_0=\sp(1)\oplus\left.
\left\{\left(\begin{array}{cc}C&0\\0&-\bar
C^t\end{array}\right)\right|C\in\Mat(r,\H)\right\}\subset\sp(1)\oplus\sp(r,r)_{W}.\end{equation}
This algebra  appears as the holonomy algebra of the
pseudo-quaternionic-K\"ahlerian symmetric space ${\rm
SL}(r+1,\H)/{\rm S}({\rm GL}(1,\H)\times {\rm GL}(r,\H))$
\cite{AC05}. This shows that $\R(\h_0)$ contains an element $R_1$
such that $\h_0$ annihilates $R_1\in\R(\h_0)$ and the image of
$R_1$ spans $\h_0$. Let $\g\subset\h_0$ and $R\in\R(\g)$.  Let
$X,Y\in W$ and $X_1\in W_1$. The Bianchi identity \eqref{Bianchi1}
implies $R(X,Y)X_1=0$ and $$R(X_1,X)Y=R(X_1,Y)X,$$ i.e.
$R(X,Y)=0$, and for each fixed $X_1\in W_1$,
$R(X_1,\cdot|_{W})|_W$ belongs to the first prolongation
$(\g|_{W})^{(1)}$ of the subalgebra
$\g|_{W}\subset\sp(1)\oplus\gl(r,\H)\subset\gl(4r,\Real)$.
Similarly, if $X_1,Y_1\in W_1$, then $R(X_1,Y_1)=0$; if $X\in W$,
then $R(X,\cdot|_{W_1})|_{W_1}\in (\g|_{W_1})^{(1)}$. It holds
$(\gl(r,\H))^{(1)}=0$ \cite{Bryant}, hence $\R(\h_0\cap
\sp(r,r))=0$. From this and  \eqref{Rsp1+}  it follows that
$\R(\h_0)=\Real R_1$. This and \eqref{Rsp1+} show that
\begin{equation}\label{Rsp1++} \R(\sp(1)\oplus\sp(r,r)_{W})=
\Real R_1\oplus\R(\sp(r,r)_{W}).\end{equation}

Let $\g\subset\sp(1)\oplus\sp(r,r)_{W }$ be a Berger subalgebra
such that $\g\not\subset\sp(r,r)_{W}$. Then there exists
$R\in\R(\g)$ such that $R=\nu R_1+\W$, $\nu\neq 0$ and
$\W\in\R(\sp(r,r)_{W})$. Let $p,X\in W$ and $Y,Z\in\Real^{4r,4r}$.
From \eqref{sym} applied to $\W$, we obtain $\W(p,X)=0$. From this
and the Bianchi identity it follows that $\W(X_1,\cdot|_{W})|_W\in
(\gl(W))^{(1)}=0$ for any $X_1\in W_1$. Hence, $R(X,X_1)|_W=\nu
R_1(X,X_1)|_W$ for all $X\in W$ and $X_1\in W_1$. This shows that
 $\g|_W=\h_0|_W=\sp(1)|_W\oplus\gl(W)$. Suppose
that $\g\neq \h_0$. Then for some $B\in S(r,\H)$ the element
$\xi=\left(\begin{array}{cc}E_r&B\\0&-E_r\end{array}\right)$
belongs to $\g$. If $B\neq 0$, then choosing the basis
$p_1,...,p_r,q'_1,...,q'_r$, where
$q'_i=q_i-\frac{1}{2}\sum_{j=1}^rB_{ji}p_j$, we get that
$\xi=\left(\begin{array}{cc}E_r&0\\0&-E_r\end{array}\right)\in\g$.
Let $\xi_1=\left(\begin{array}{cc}C&B\\0&-\bar
C^t\end{array}\right)\in\g$, where  $C\in\Mat(r,\H)$ and $B\in
S(r,\H)$. Then,
$[\xi,\xi_1]=\left(\begin{array}{cc}0&2B\\0&0\end{array}\right)\in\g$.
This shows that $\g=\h_0\zr L$, where $L\subset\left.
\left\{\left(\begin{array}{cc}0&B\\0&0\end{array}\right)\right|B\in
S(r,\H)\right\}$. The Lie brackets of elements from $\h_0$ and $L$
are given by the representation of $\gl(r,\H)$ on $S(r,\H)$. Since
this representation is irreducible \cite{Bryant},
$\g=\sp(1)\oplus\sp(r,r)_W$. It is not hard to see that this
algebra is a Berger algebra, hence it is a candidate to be a
holonomy algebra. $\Box$

\section{Proof of Proposition \ref{prop1}}
Let $(M,g)$ be  pseudo-quaternionic-K\"ahlerian manifold with the
holonomy algebra  $\h_0$. Fix a point $x\in M$. Then $T_xM$ is
identified with $\Real^{4r,4r}$. For the covariant derivative of
the curvature tensor at the point $x$ we have $\nabla_X
R_x\in\R(\h_0)$ for any $X\in\Real^{4r,4r}$. Let $X,Y\in W$ and
$X_1\in W_1$. From the above, we get $\nabla_{X_1}R_x(X,Y)=0$.
This and the second Bianchi identity imply
$$\nabla_{X}R_x(X_1,Y)=\nabla_{Y}R_x(X_1,X),$$ i.e.
$\nabla_{\cdot|_W}R_x(X_1,\cdot|_W)|_W$ belongs to the second
prolongation of the subalgebra
$\sp(1)\oplus\gl(r,\H)\subset\gl(4r,\Real),$ which is trivial,
since the second prolongation of its complexification
$\sl(2,\Co)\oplus\gl(2r,\Co)\subset\gl(4r,\Co)$ is trivial
\cite{Sch}. We get that $\nabla_{X}R_x=0$. By the same arguments,
$\nabla_{X_1}R_x=0$. Thus, $\nabla R_x=0$ for any $x\in M$, i.e.
$\nabla R=0$ and $(M,g)$ is locally symmetric. $\Box$

\section{Proof of corollaries}

Corollary 1 follows from the above results and from the fact that
the only irreducible holonomy group of not locally symmetric
pseudo-quaternionic-K\"ahlerian manifolds is $\Sp(1)\cdot\Sp(r,s)$
\cite{Bryant}. Corollary 3 follows from Corollary 1  and from the
results of \cite{AC05}. To prove Corollary 2 consider the
universal  covering $(\tilde M,\tilde g)$ of $(M,g)$. Then
$(\tilde M,\tilde g)$ is simply connected and complete, hence it
is contained in the list of symmetric spaces from \cite{AC05}
given in Introduction. It is known that $M=\tilde M/\Gamma$ for
some freely acting discrete group $\Gamma$ of isometries of
$\tilde M$, e.g. \cite{Borel}. $\Box$

\begin{acknowledgements}  I am grateful to D.~V.~Alekseevsky and Jan
Slov\'ak for useful discussions,  support and help. The author has
been supported by the grant GACR 201/09/H012.
\end{acknowledgements}


\begin{thebibliography}{}

\bibitem{AC05} D.~V.~Alekseevsky, V.~Cort\'es, {\it Classification of
pseudo-Riemannian symmetric spaces of quaternionic K\"ahler type.}
Lie groups and invariant theory, 33--62,
 Amer. Math. Soc. Transl. Ser. 2, 213, Amer. Math. Soc., Providence, RI, 2005.


\bibitem{Ber} M.~Berger, {\it Sur les groupers d'holonomie des
vari\'et\'es \`aconnexion affine et des vari\'et\'es
riemanniennes}. Bull. Soc. Math. France 83 (1955), 279--330.


\bibitem{Besse} A.~L.~Besse, {\it Einstein manifolds}, Springer-Verlag,
Berlin-Heidelberg-New York, 1987.

\bibitem{Borel} A.~Borel, {\it Semisimple groups and Riemannian symmetric spaces.}
 Texts and Readings in Mathematics, 16. Hindustan Book Agency, New Delhi, 1998. x+136 pp.

\bibitem{Bryant} R.~L.~Bryant, {\it Classical, exceptional, and exotic
holonomies: a status report.} Actes de la Table Ronde
 de G\'eom\'etrie Diff\'erentielle
 (Luminy, 1992), 93--165, Semin. Congr., 1, Soc. Math. France, Paris, 1996.


\bibitem{ESI} A.~S.~Galaev, T.~Leistner, {\it Holonomy groups of Lorentzian manifolds: classification,
examples, and applications.} Recent developments in
pseudo-Riemannian geometry, 53--96, ESI Lect. Math. Phys., Eur.
Math. Soc., Z\"urich, 2008.

\bibitem{IRMA} A.~S.~Galaev, T.~Leistner, {\it
Recent developments in pseudo-Riemannian holonomy theory.}
 Cort\'es, Vicente (ed.), Handbook of pseudo-Riemannian
geometry and supersymmetry.  Zurich: European Mathematical
Society. IRMA Lectures in Mathematics and Theoretical Physics 16,
581-627 (2010).

\bibitem{Jo} D.~Joyce, {\it Riemannian holonomy  groups and calibrated
geometry}. Oxford University Press, 2007.


\bibitem{March} {\it Quaternionic structures in mathematics and physics.} Proceedings of
the 2nd Meeting held in Rome, September 6--10, 1999. Ed.
S.~Marchiafava, P.~Piccinni and M.~Pontecorvo. World Scientific
Publishing Co., Inc., River Edge, NJ, 2001. xvi+469 pp.


 \bibitem{Salamon} S.~M.~Salamon, {\it Differential geometry of quaternionic manifolds.}
  Ann. Sci. \'Ecole Norm. Sup. (4) 19 (1986), no. 1, 31--55.

\bibitem{Sch} L.~J.~Schwachh\"ofer, {\it Connections with irreducible holonomy representations,}
 Adv. Math. 160 (2001), no. 1, 1--80.


\bibitem{Swann} A.~Swann, {\it Hyper-K\"ahler and quaternionic K\"ahler geometry.} Math. Ann. 289
(1991), no. 3, 421--450.



\end{thebibliography}


\end{document}